\documentclass[11pt]{amsart}
\usepackage{amssymb}
\usepackage[margin=1in]{geometry}
\usepackage[colorlinks]{hyperref}
\newtheorem{theorem}{Theorem}[section]

\theoremstyle{definition}

\theoremstyle{remark}
\newtheorem{remark}[theorem]{Remark}

\numberwithin{equation}{section}

\begin{document}

\newcommand{\spacing}[1]{\renewcommand{\baselinestretch}{#1}\large\normalsize}
\spacing{1.14}
\title[On the classification of five-dimensional nilsolitons]{On the classification of five-dimensional nilsolitons}

\author {Hamid Reza Salimi Moghaddam}

\address{Department of Pure Mathematics \\ Faculty of  Mathematics and Statistics\\ University of Isfahan\\ Isfahan\\ 81746-73441-Iran.\\ Scopus Author ID: 26534920800 \\ ORCID Id:0000-0001-6112-4259\\} \email{hr.salimi@sci.ui.ac.ir and salimi.moghaddam@gmail.com}


%%\date{\today}

\begin{abstract}
In 2002, using a variational method, Lauret classified five-dimensional nilsolitons. In this work, using the algebraic Ricci soliton equation, we obtain the same classification. We show that, among ten classes of five-dimensional nilmanifolds, seven classes admit Ricci soliton structure. In any case, the derivation which satisfies the algebraic Ricci soliton equation is computed.
\end{abstract}

\maketitle

%%-------------------------Introduction-------------------------
\section{\textbf{Introduction}}\label{Introduction}
Suppose that $(M,g)$ is a complete Riemannian manifold and $\textsf{ric}_g$ denotes its Ricci tensor. If, for a real number $c$ and a complete vector field $X$, the Riemannian metric $g$ satisfies the equation
\begin{equation}\label{Main Ricci soliton equation}
    \textsf{ric}_g=c g+\textsf{L}_Xg,
\end{equation}
then $g$ is called a Ricci soliton. In this equation, if $c>0$, $c=0$ or $c<0$, then $g$ is named shrinking, steady, or expanding Ricci soliton. \\
Although the above definition is a generalization of Einstein metrics but the main motivation for considering Ricci solitons is the Ricci flow equation
\begin{equation}\label{Ricci flow equation}
    \frac{\partial}{\partial t}g_t=-2\textsf{ric}_{g(t)},
\end{equation}
where $\phi_t$ is a one-parameter group of diffeomorphisms. A Riemannian metric $g$ is a Ricci soliton on $M$ if and only if the following one-parameter family of Riemannian metrics is a solution of (\ref{Ricci flow equation}) (see \cite{Lauret1} and \cite{Lauret2}),
\begin{equation}\label{family of Riemannian metrics}
    g_t=(-2c t+1)\phi_t^\ast g.
\end{equation}
A very interesting case happens when we study Ricci solitons on Lie groups. Suppose that $G$ is a Lie group, $\frak{g}$ is its Lie algebra and $g$ is a left invariant Riemannian metric on $G$. If for a real number $c$ and a derivation $D\in\textsf{Der}(\frak{g})$, the $(1,1)$-Ricci tensor $\textsf{Ric}_g$ of $g$ satisfies the equation
\begin{equation}\label{Algebraic Ricci soliton equation}
    \textsf{Ric}_g=c\textsl{Id}+D,
\end{equation}
then $g$ is called an algebraic Ricci soliton. Moreover, if $G$ is a nilpotent Lie group then it is called a nilsoliton. It is shown that all algebraic Ricci solitons are Ricci soliton (for more details see \cite{Jablonski1} and \cite{Lauret3}).\\
In the year 2001, Lauret proved that on nilpotent Lie groups a left invariant Riemannian metric $g$ is a Ricci soliton if and only if it is an algebraic Ricci soliton. In fact for left invariant Riemannian metrics on nilpotent Lie groups the equations (\ref{Main Ricci soliton equation}) and (\ref{Algebraic Ricci soliton equation}) are equivalent (see \cite{Lauret3} and \cite{Jablonski2}).\\
Naturally, one can generalize the concept of algebraic Ricci soliton to homogeneous spaces. In 2014, Jablonski showed that any homogeneous Ricci soliton is an algebraic Ricci soliton (see \cite{Jablonski2}). So the classification of left invariant Ricci solitons on Lie groups reduces to the classification of algebraic Ricci solitons on them. \\
In the year 2002, using a variational method, Lauret classified five-dimensional nilsolitons (for more details see \cite{Lauret4} and \cite{Will}). In this work, we classify five-dimensional nilsolitons using the classification of five-dimensional nilmanifolds given in \cite{Homolya-Kowalski} and \cite{Figula-Nagy}, and the algebraic Ricci soliton equation (\ref{Algebraic Ricci soliton equation}). We see that our results are compatible with the results of \cite{Lauret4}.\\
Suppose that $g$ is a left invariant Riemannian metric on a Lie group $G$ and $\alpha_{ijk}$ are the structural constants of the Lie algebra $\frak{g}$, with respect to an orthonormal basis $\{E_1,\cdots,E_n\}$, defined by the following equations:
\begin{equation}\label{structural constants}
    [E_i,E_j]=\sum_{k=1}^n\alpha_{ijk}E_k.
\end{equation}
In an earlier paper (see \cite{Salimi}), we proved that $g$ is an algebraic Ricci soliton if and only if there exists a real number $c$ such that, for any $t, p, q=1,\cdots,n$ the structural constants satisfie in the following equation:
\begin{eqnarray}\label{Main formula}
  c\alpha_{qpt}+\frac{1}{4} \sum_{i=1}^n\sum_{j=1}^n\sum_{r=1}^n &&\hspace{-0.6cm} 2\alpha_{rjj}\Big{(}\alpha_{iqt}(\alpha_{pri}+\alpha_{ipr}-\alpha_{rip})
  -\alpha_{ipt}(\alpha_{qri}+\alpha_{iqr}-\alpha_{riq})\Big{)}\nonumber\\
  &&\hspace{-0.6cm} +2(\alpha_{rji}+\alpha_{irj}-\alpha_{jir})(\alpha_{ipt}\alpha_{qjr}-\alpha_{iqt}\alpha_{pjr})\nonumber\\
  &&\hspace{-0.6cm} +(\alpha_{jri}+\alpha_{ijr}-\alpha_{rij})\Big{(}\alpha_{ipt}(\alpha_{qjr}+\alpha_{rqj}-\alpha_{jrq})
  -\alpha_{iqt}(\alpha_{pjr}+\alpha_{rpj}-\alpha_{jrp})\Big{)}\\
  &&\hspace{-0.6cm} -2\alpha_{pqi}\alpha_{rjj}(\alpha_{irt}+\alpha_{tir}-\alpha_{rti})+2\alpha_{pqi}\alpha_{ijr}(\alpha_{rjt}+\alpha_{trj}-\alpha_{jtr})\nonumber\\
  &&\hspace{-0.6cm} +\alpha_{pqi}(\alpha_{ijr}+\alpha_{rij}-\alpha_{jri})(\alpha_{jrt}+\alpha_{tjr}-\alpha_{rtj})=0. \nonumber
\end{eqnarray}
In this case we showed that, for any $i=1,\cdots,n$, the derivation $D$ satisfies in the following equation:
\begin{eqnarray}\label{Derivation equation}
  D(E_i) &=& -cE_i+\frac{1}{4} \sum_{l=1}^n\Big{\{} \sum_{j=1}^n\sum_{r=1}^n2\alpha_{rjj}(\alpha_{irl}+\alpha_{lir}-\alpha_{rli})\nonumber \\
         &&\hspace{4cm}-2\alpha_{ijr}(\alpha_{ril}+\alpha_{lrj}-\alpha_{jlr})\\
         &&\hspace{4cm}-(\alpha_{ijr}+\alpha_{rij}-\alpha_{jri})(\alpha_{jrl}+\alpha_{ljr}-\alpha_{rlj})\Big{\}}E_l.\nonumber
\end{eqnarray}
In \cite{Salimi}, using the above equations, we classified all algebraic Ricci solitons on three-dimensional Lie groups. \\
Using the equations (\ref{Main formula}) and (\ref{Derivation equation}) which is different from Lauret's variational method, again, we classify all five-dimensional nilsolitons. We remember that in this method, against some recent works (such as \cite{Di Cerbo}) which have studied the Ricci solitons only on the set of left invariant vector fields, we consider the Ricci solitons on all vector fields.\\
A Riemannian nilmanifold is a connected Riemannian manifold $(M,g)$ such that there exists a nilpotent Lie subgroup of its isometry group $I(M)$ which acts transitively on $M$. Wilson, in \cite{Wilson}, proved that if $M$ is a homogeneous Riemannian nilmanifold then there is a unique nilpotent Lie subgroup of $I(M)$ which acts simply transitively on $M$. Also, he showed that this Lie subgroup is a normal subgroup of $I(M)$. Thus, we can assume a homogeneous Riemannian nilmanifold as a nilpotent Lie group $N$ equipped with a left invariant Riemannian metric $g$. \\
In 2006, Homoloya and Kowalski, up to isometry, classified five-dimensional two-step nilpotent Lie groups equipped with left invariant Riemannian metrics (see \cite{Homolya-Kowalski}). In a paper written by Figula and Nagy in 2018, up to isometry, the classification of five-dimensional nilmanifolds has been completed by the classification of five-dimensional nilpotent Lie algebras of nilpotency classes three and four equipped with inner products. We mention that this classification does not contain the Lie algebras which are direct products of Lie algebras of lower dimensions (see \cite{Figula-Nagy}).\\
In the next section we will classify all left invariant Ricci solitons on all ten classes of five-dimensional nilmanifolds.

%%--------The classification of algebraic Ricci solitons on five-dimensional nilmanifolds-----------
\section{\textbf{The classification of five-dimensional nilsolitons}}\label{The classification}
In \cite{Homolya-Kowalski}, all five-dimensional two-step nilmanifolds up to isometry have been classified  by Homolya and Kowalski. Recently, in \cite{Figula-Nagy}, Figula and Nagy have completed the classification of five-dimensional nilmanifolds up to isometry by the classification of five-dimensional nilpotent Lie groups of nilpotency class greater than two. In this section using formula (\ref{Main formula}) and the above classifications we classify all algebraic Ricci solitons on simply connected five-dimensional nilmanifolds. \\
Suppose that $N$ is an arbitrary simply connected five-dimensional nilmanifold then, up to isometry, its Lie algebra is one of the following ten cases. In all cases the set $\{E_1, \cdots, E_5\}$ is an orthonormal basis for the Lie algebra and in any case we give only non-vanishing commutators. We use the notation $d(a_1,\cdots,a_n)$ to denote the diagonal matrix with entries $a_1,\cdots,a_n$.\\
In the cases of nilpotency class greater than two, we will consider in detail the cases 2.6 and 2.7, which are the most difficult cases.

\subsection{\textbf{Two-step nilpotent Lie algebra with one-dimensional center.}}
Let $N$ be the two-step nilpotent Lie group with one-dimensional center and $\frak{n}$ denotes its Lie algebra. By \cite{Homolya-Kowalski}, the non-zero Lie brackets are as follows:
\begin{equation}
    [E_1,E_2]=sE_5, \ \ \ \ [E_3,E_4]=mE_5,
\end{equation}
where $s\geq m >0$. So, in this case, the structure constants with respect to the orthonormal basis $\{E_1, \cdots, E_5\}$ are of the following forms:
\begin{equation*}
    \alpha_{125}=-\alpha_{215}=s, \ \ \ \ \alpha_{345}=-\alpha_{435}=m.
\end{equation*}
Easily, using the formula (\ref{Main formula}), we see that $\frak{n}$ is an algebraic Ricci soliton if and only if
\begin{equation*}
    \left\{
  \begin{array}{ll}
    cs+\frac{1}{2}sm^2+\frac{3}{2}s^3=0, &  \\
    cm+\frac{1}{2}ms^2+\frac{3}{2}m^3=0. &
  \end{array}
\right.
\end{equation*}
A direct computation shows that the above equations hold if and only if $m=s$. So for the constant $c$ we have $c=-2m^2$. Now, the equation (\ref{Derivation equation}) shows that for the matrix representation of $D$ we have $$D=d(\frac{3}{2}m^2,\frac{3}{2}m^2,\frac{3}{2}m^2,\frac{3}{2}m^2,3m^2).$$

\subsection{\textbf{Two-step nilpotent Lie algebra with two-dimensional center.}}
The second case is two-step nilpotent Lie group $N$ with two-dimensional center. In \cite{Homolya-Kowalski}, it is shown that there are real numbers $m\geq s >0$ such that the non-zero Lie brackets are
\begin{equation}
    [E_1,E_2]=mE_4, \ \ \ \ [E_1,E_3]=sE_5.
\end{equation}
In this case for the structure constants we have:
\begin{equation*}
    \alpha_{124}=-\alpha_{214}=m, \ \ \ \ \alpha_{135}=-\alpha_{315}=s.
\end{equation*}
Now, the equation (\ref{Main formula}) shows that $\frak{n}$ satisfies the algebraic Ricci soliton equation (\ref{Algebraic Ricci soliton equation}) if and only if
\begin{equation*}
    \left\{
  \begin{array}{ll}
    cm+\frac{3}{2}m^3+\frac{1}{2}s^2m=0, &  \\
    cs+\frac{3}{2}s^3+\frac{1}{2}m^2s=0. &
  \end{array}
\right.
\end{equation*}
Easily we can see, for the solutions we have $m=s$ and $c=-2m^2$. Now, by equation (\ref{Derivation equation}), for the matrix representation of $D$ in the basis $\{E_1, \cdots, E_5\}$ we have $$D=d(m^2,\frac{3}{2}m^2,\frac{3}{2}m^2,\frac{5}{2}m^2,\frac{5}{2}m^2).$$

\subsection{\textbf{Two-step nilpotent Lie algebra with three-dimensional center.}}
The third Lie group which we have considered is the two-step nilpotent Lie group $N$ with three-dimensional center. This is the last five-dimensional two-step nilpotent Lie group.  Based on \cite{Homolya-Kowalski}, the non-zero Lie bracket of this case is
\begin{equation}
    [E_1,E_2]=mE_3,
\end{equation}
where $m$ is a positive real number. So, the non-zero structure constants are
\begin{equation*}
    \alpha_{123}=-\alpha_{213}=m.
\end{equation*}
It is easy to show that the equation (\ref{Main formula}) holds if and only if $cm+\frac{2}{3}m^3=0$. The last equation shows that $c=-\frac{3}{2}m^2$. So, by (\ref{Derivation equation}), the representation of $D$ in the basis $\{E_1, \cdots, E_5\}$, is of the form $$D=d(m^2,m^2,2m^2,\frac{3}{2}m^2,\frac{3}{2}m^2).$$

\subsection{\textbf{Four-step nilpotent Lie algebra $\frak{l}_{5,7}$, (case A)}}
In \cite{Figula-Nagy}, Figula and Nagy have proven that the five-dimensional nilpotent Lie algebras $\frak{n}$ of nilpotency class greater than two are of the forms $\frak{l}_{5,7}$, $\frak{l}_{5,6}$, $\frak{l}_{5,5}$ and $\frak{l}_{5,9}$. In this subsection and the next subsection, we study the necessary and sufficient conditions for the Lie algebra $\frak{l}_{5,7}$ to be an algebraic Ricci soliton. By \cite{Figula-Nagy}, the non-vanishing Lie brackets of $\frak{l}_{5,7}$ are
\begin{equation}
    [E_1,E_2]=mE_3+sE_4+uE_5, \ \ \ [E_1,E_3]=vE_4+wE_5, \ \ \ [E_1,E_4]=xE_5,
\end{equation}
where for the real numbers $m,s,u,v,w,x$ we have two cases: (case A: $m,v,x>0$, $s=0$ and $w\geq 0$) and (case B: $m,v,x>0$ and $s>0$). In this subsection we study the case A.
In the case of $\frak{l}_{5,7}$ (in two cases A and B) for the structure constants we have:
\begin{eqnarray*}
  &&\alpha_{123}=-\alpha_{213}=m, \ \alpha_{124}=-\alpha_{214}=s, \ \alpha_{125}=-\alpha_{215}=u, \\
  &&\alpha_{134}=-\alpha_{314}=v, \ \alpha_{135}=-\alpha_{315}=w, \ \alpha_{145}=-\alpha_{415}=x.
\end{eqnarray*}

The equation (\ref{Main formula}) together with some computations  shows that $\frak{l}_{5,7} (case A)$, satisfies the algebraic Ricci soliton equation (\ref{Algebraic Ricci soliton equation}) if and only if
\begin{equation*}
    x=m, \ \ c=-2m^2, \ \ u=w=s=0, \  and \ \ v=\frac{2}{\sqrt{3}}m.
\end{equation*}

Now, using (\ref{Derivation equation}), the matrix representation of $D$ is of the form $$D=d(\frac{1}{3}m^2,\frac{3}{2}m^2,\frac{11}{6}m^2,\frac{13}{6}m^2,\frac{5}{2}m^2).$$

\subsection{\textbf{Four-step nilpotent Lie algebra $\frak{l}_{5,7}$, (case B)}}
As we mentioned above, in this case (case B) for the same structure constants of case A, we have $m,v,x>0$ and $s>0$. We see that the equation (\ref{Main formula}) deduces to a system of equations with no solution. So the Lie algebra $\frak{l}_{5,7}$, (case B) does not admit algebraic Ricci soliton structure.

\subsection{\textbf{Four-step nilpotent Lie algebra $\frak{l}_{5,6}$, (case A)}}
During this subsection and the next subsection, we have considered the Lie algebra $\frak{l}_{5,6}$. In two cases (cases A and B), by \cite{Figula-Nagy}, the Lie brackets are
\begin{equation}
    [E_1,E_2]=mE_3+sE_4+uE_5, \ \ \ [E_1,E_3]=vE_4+wE_5, \ \ \ [E_1,E_4]=xE_5, \ \ \ [E_2,E_3]=yE_5,
\end{equation}
where $m,s,u,v,w,x,y$ are real numbers such that $m,v,x,y\neq0$.
For the structure constants we have:
\begin{eqnarray*}
  &&\alpha_{123}=-\alpha_{213}=m, \ \alpha_{124}=-\alpha_{214}=s, \ \alpha_{125}=-\alpha_{215}=u, \ \alpha_{134}=-\alpha_{314}=v,\\
  &&\alpha_{135}=-\alpha_{315}=w, \ \alpha_{145}=-\alpha_{415}=x, \ \alpha_{235}=-\alpha_{325}=y.
\end{eqnarray*}
Similar to the Lie algebra $\frak{l}_{5,7}$ we have two cases, case A and B. In the case A we have $s=0$ and $w\geq 0$ and for the case B $s>0$.

A direct computation, using (\ref{Main formula}), shows that the Lie algebra $\frak{l}_{5,6} (case A)$, is an algebraic Ricci soliton if and only if (for two cases A and B)
\begin{equation*}
    \left\{
  \begin{array}{ll}
    myu=0, &  \\
    mvs+mwu+sxu=0, &  \\
    cm+\frac{3m}{2}(m^2+s^2+u^2)+\frac{x}{2}(xm-sw)=0, &  \\
    cs+\frac{3s}{2}(m^2+s^2+u^2+v^2)+\frac{1}{2}(-mwx+w^2s+y^2s)+wvu=0, &  \\
    cu+\frac{3u}{2}(m^2+s^2+u^2+w^2+x^2+y^2)+\frac{v^2u}{2}+vsw=0, &  \\
    uxv=0, &  \\
    mvs+mwu-\frac{vxw}{2}=0, &  \\
    cv+\frac{3v}{2}(s^2+v^2+w^2)+\frac{v}{2}(u^2+y^2)+usw=0, &  \\
    cw+\frac{3w}{2}(u^2+v^2+w^2+x^2+y^2)+\frac{s}{2}(-mx+sw)+svu=0, & \\
    umx=0, &  \\
    sxu+vxw-\frac{mvs}{2}=0, & \\
    cx+\frac{3x}{2}(u^2+w^2+x^2)+\frac{1}{2}(m^2x-msw+xy^2)=0, &  \\
    mvu=0, & \\
    mwu+sxu+vxw=0, &  \\
    suy+vwy=0, &  \\
    cy+\frac{3y}{2}(u^2+w^2+y^2)+\frac{y}{2}(s^2+v^2+x^2)=0, &  \\
    wyx-\frac{msy}{2}=0, &  \\
    uyx=0. &
  \end{array}
\right.
\end{equation*}
A direct computation shows that the above equations hold if and only if
\begin{equation*}
    u=w=s=0, m=\pm\sqrt{\frac{3}{2}}x, y=\pm x, v=\pm\sqrt{\frac{3}{2}}x \ and \ c=-\frac{11}{4}x^2.
\end{equation*}

The equation (\ref{Derivation equation}) shows that for the matrix representation of $D$ we have $$D=d(\frac{3}{4}x^2,\frac{3}{2}x^2,\frac{9}{4}x^2,3x^2,\frac{15}{4}x^2).$$

\subsection{\textbf{Four-step nilpotent Lie algebra $\frak{l}_{5,6}$, (case B)}}
In the case B for the same structure constants of the above case, we have $s>0$. Then we see that the system of equations defined by (\ref{Main formula}) has no solution. Therefore similar to the case $\frak{l}_{5,7}$ (case B), the Lie algebra $\frak{l}_{5,6}$ (case B), does not admit an algebraic Ricci soliton structure.

\subsection{\textbf{Three-step nilpotent Lie algebra $\frak{l}_{5,5}$}}
This section is devoted to $\frak{l}_{5,5}$ which is a three-step nilpotent Lie algebra with one dimensional center.
For the real numbers $s,u\geq0$, and $m,v,w>0$, the non-vanishing Lie brackets are as follows:
\begin{equation}
    [E_1,E_2]=mE_4+sE_5, \ \ \ [E_1,E_3]=uE_5, \ \ \ [E_1,E_4]=vE_5, \ \ \ [E_2,E_3]=wE_5.
\end{equation}
Hence for the structure constants we have:
\begin{eqnarray*}
  &&\alpha_{124}=-\alpha_{214}=m, \ \alpha_{125}=-\alpha_{215}=s, \ \alpha_{135}=-\alpha_{315}=u, \\ &&\alpha_{145}=-\alpha_{145}=v, \ \alpha_{235}=-\alpha_{325}=w.
\end{eqnarray*}
Now substituting the above structure constants in the equation (\ref{Main formula}), leads us to a system of equations with the following solution,
\begin{equation*}
    u=s=0, v=m, w=\frac{\sqrt{2}}{2}m, c=-\frac{7}{4}m^2.
\end{equation*}
So the Lie algebra $\frak{l}_{5,5}$ admits an algebraic Ricci soliton structure if and only if the above equations hold. In this case, by (\ref{Derivation equation}), the matrix representation of $D$ in the basis $\{E_1, \cdots, E_5\}$, is $$D=d(\frac{3}{4}m^2,m^2,\frac{3}{2}m^2,\frac{7}{4}m^2,\frac{5}{2}m^2).$$

\subsection{\textbf{Three-step nilpotent Lie algebra $\frak{l}_{5,9}$, (case A)}}
For the non-vanishing Lie brackets of this case we have:
\begin{equation}
    [E_1,E_2]=mE_3+sE_4+uE_5, \ \ \ [E_1,E_3]=vE_4, \ \ \ [E_2,E_3]=wE_5,
\end{equation}
where $m>0$, $w>v>0$ and $s,u\geq0$.
So we have:
\begin{eqnarray*}
  &&\alpha_{123}=-\alpha_{213}=m, \ \alpha_{124}=-\alpha_{214}=s, \ \alpha_{125}=-\alpha_{215}=u, \\ &&\alpha_{134}=-\alpha_{314}=v, \ \alpha_{235}=-\alpha_{325}=w.
\end{eqnarray*}
Easily we see that the system of equations which is defied by (\ref{Main formula}) has no solution. Therefore, the case A of the Lie algebra $\frak{l}_{5,9}$ does not admit Algebraic Ricci soliton structure.

\subsection{\textbf{Three-step nilpotent Lie algebra $\frak{l}_{5,9}$, (case B)}}
In the case B of the Lie algebra $\frak{l}_{5,9}$, for the same Lie algebra brackets, we have $m,v>0$, $v=w$, $u=0$ and $s\geq0$. Then, the equation (\ref{Main formula}), leads us to the following system of equations:
\begin{equation*}
    \left\{
  \begin{array}{l}
    mvs=0,\\
    cm+\frac{3}{2}(m^3+s^2m)=0,  \\
    cs+\frac{3}{2}(s^3+m^2s)+2v^2s=0,\\
    cv+\frac{3}{2}s^2v+2v^3=0,\\
    cv+\frac{1}{2}s^2v+2v^3=0.
  \end{array}
\right.
\end{equation*}
The solution of the above system of equations is as follows:
\begin{equation*}
    s=0, v=\frac{\sqrt{3}}{2}m, c=-\frac{3}{2}m^2.
\end{equation*}
So the case B of the Lie algebra $\frak{l}_{5,9}$ admits an algebraic Ricci soliton structure with the derivation $$D=d(\frac{5}{8}m^2,\frac{5}{8}m^2,\frac{5}{4}m^2,\frac{15}{8}m^2,\frac{15}{8}m^2).$$

Now, we summarize the above results and give their relations with Lauret classification (Theorem 5.1 of \cite{Lauret4}) in the table \ref{The summarized result}.

{\fontsize{7}{0}{\selectfont

\begin{table}
    \centering\caption{The classification of five-dimensional nilsolitons}\label{The summarized result}
       \begin{tabular}{|p{0.5cm}|p{3.5cm}|p{1.5cm}|p{1.5cm}|p{5cm}|p{1.5cm}|}
        \hline
            case & Lie algebra structure & conditions for parameters of  Lie algebra structure \newline & Ricci soliton & the constant $c$ and the derivation $D$ & equivalence class in Lauret classification \\
             \hline
             case 2.1 & $[E_1,E_2]=sE_3$, \newline $[E_3,E_4]=mE_5$ & $s\geq m>0$ & $+$, if $s=m$ & $c=-2m^2$,\newline $D=d(\frac{3}{2}m^2,\frac{3}{2}m^2,\frac{3}{2}m^2,\frac{3}{2}m^2,3m^2)$ & $\mu_4'$ \\
            \hline
            case 2.2 & $[E_1,E_2]=mE_4$, \newline $[E_1,E_3]=sE_5$ & $m\geq s>0$ & $+$, if $s=m$ & $c=-2m^2$,\newline $D=d(m^2,\frac{3}{2}m^2,\frac{3}{2}m^2,\frac{5}{2}m^2,\frac{5}{2}m^2)$ & $\mu_6'$ \\
            \hline
            case 2.3 & $[E_1,E_2]=mE_3$ & $m>0$ & $+$ & $c=-\frac{3}{2}m^2$,\newline $D=d(m^2,m^2,2m^2,\frac{3}{2}m^2,\frac{3}{2}m^2)$ & $\mu_7'$ \\
            \hline
            case 2.4 & $[E_1,E_2]=mE_3+sE_4+uE_5$,\newline  $[E_1,E_3]=vE_4+wE_5$, \newline $[E_1,E_4]=xE_5$ & $m,v,x>0$, $s=0$ and $w\geq0$ & $+$, if $x=m$, $u=w=s=0$ and $v=\frac{2}{\sqrt{3}}m$ & $c=-2m^2$,\newline $D=d(\frac{1}{3}m^2,\frac{3}{2}m^2,\frac{11}{6}m^2,\frac{13}{6}m^2,\frac{5}{2}m^2)$ & $\mu_1'$ \\
            \hline
            case 2.5 & $"$ & $m,v,x>0$, $w\geq0$ and $s>0$ & $-$ & $-$ & $-$ \\
            \hline
            case 2.6 & $[E_1,E_2]=mE_3+sE_4+uE_5$,\newline  $[E_1,E_3]=vE_4+wE_5$, \newline $[E_1,E_4]=xE_5$, \newline $[E_2,E_3]=yE_5$ & $m,v,x,y\neq0$,\newline $s=0$, $w\geq 0$ & $+$, if $u=w=s=0$ \newline $m=v=\pm\sqrt{\frac{3}{2}}x$, \newline $y=\pm x$ & $c=-\frac{11}{4}x^2$,\newline $D=d(\frac{3}{4}x^2,\frac{3}{2}x^2,\frac{9}{4}x^2,3x^2,\frac{15}{4}x^2)$ & $\mu_2'$ \\
            \hline
            case 2.7 & $"$ & $s>0$ & $-$ & $-$ & $-$ \\
            \hline
            case 2.8 & $[E_1,E_2]=mE_4+sE_5$,\newline  $[E_1,E_3]=uE_5$, \newline $[E_1,E_4]=vE_5$, \newline $[E_2,E_3]=wE_5$ & $s,u\geq0$, $m,v,w>0$ & $+$, if $s=u=0$, $v=m$ and $w=\frac{\sqrt{2}}{2}m$ & $c=-\frac{7}{4}m^2$,\newline $D=d(\frac{3}{4}m^2,m^2,\frac{3}{2}m^2,\frac{7}{4}m^2,\frac{5}{2}m^2)$ & $\mu_3'$ \\
            \hline
            case 2.9 & $[E_1,E_2]=mE_3+sE_4+uE_5$,\newline  $[E_1,E_3]=vE_4$, \newline $[E_2,E_3]=wE_5$, & $m>0$, $w>v>0$, $s,u\geq0$ & $-$ & $-$ & $-$ \\
            \hline
            case 2.10 & $"$ & $m,v>0$,\newline $w=v$, \newline $s\geq0$,\newline $u=0$ & $+$, if $s=0$,\newline $v=\frac{\sqrt{3}}{2}m$ & $c=-\frac{3}{2}m^2$,\newline $D=d(\frac{5}{8}m^2,\frac{5}{8}m^2,\frac{5}{4}m^2,\frac{15}{8}m^2,\frac{15}{8}m^2)$ & $\mu_5'$ \\
            \hline
          \end{tabular}
        \end{table}
}}
\begin{remark}
For the Lie groups of nilpotency class greater than two we have used the classification given in \cite{Figula-Nagy}, and we only consider Lie algebras which are not direct products of Lie algebras of lower dimension. Therefore, class $\mu_8'$ of Lauret classification does not appear in our classification.
\end{remark}
\begin{remark}
Based on the results given in table \ref{The summarized result}, there are some three and four-step nilmanifolds which are not nilsoliton.
\end{remark}
%%-------------------- BIBLIOGRAPHY------------------------

\bibliographystyle{amsplain}

\end{document}